# There is no non-zero Stable Fixed Point for dense networks in the homogeneous Kuramoto model

Richard Taylor

Australian Defence Science and Technology Organisation

E-mail: richard.taylor@dsto.defence.gov.au

**Abstract**

This paper is concerned with the existence of multiple stable fixed point solutions of the homogeneous Kuramoto model. We develop a necessary condition for the existence of stable fixed points for the general network Kuramoto model. This condition is applied to show that for sufficiently dense $n$-node networks, with node degrees at least $0.9395(n-1)$, the homogeneous (equal frequencies) model has only one stable fixed point solution over the full space of phase angles in the range $-\pi$ to $\pi$. This is the zero fixed point solution defined by all phase angle differences being zero. This result together with existing research proves a conjecture of Verwoerd and Mason (2007) that for the complete network and homogeneous model the zero fixed point has a basin of attraction consisting of the entire space minus a set of measure zero. The necessary conditions are also tested to see how close to sufficiency they might be by applying them to a class of regular degree networks studied by Wiley, Strogatz and Girvan (2006).

## 1 Introduction

The Kuramoto model (see Kuramoto (1975)) was originally motivated by the phenomenon of collective synchronization whereby a system of oscillating entities (or nodes) will sometimes synchronize despite differences in the natural frequencies of the individual nodes. The model is relevant to a number of phenomena including ecology, biology, physics and social and organizational systems. Strogatz (2000) provides an accessible introduction and a number of surveys summarise many technical results (see Acebron *et al* (2005), Arenas *et al* (2008), Boccaletti *et al* (2006), Dorogovstev *et al* (2008)). The results about the Kuramoto model focus on a number of areas. These include the conditions under which a network can (partially or completely) synchronize, the behaviour of the system near synchronization, and the geometry (size and position) of attractor regions to which the initial states of a network can converge.

Kuramoto (1975) studied the infinite complete network (in which all nodes are linked) and was able to derive elegant results about the synchronization of such networks. Since the work of Kuramoto, many researchers have also studied finite networks that might better reflect real world situations. Among finite networks the complete networks on a given number of nodes form a natural finite counterpart to the infinite complete networks studied by Kuramoto. Finite networks with a general topology have also been studied with a view to understanding the role of the topology in the dynamics and synchronization of the network. Topology would model particular node relationships between a subset of all node pairs, consistent with social or organizational relationships for example. This paper is confined to the study of the Kuramoto model in the context of finite networks, with particular results for dense networks.

*1.1 The model*

The basic governing equation is given by

$$\dot{\theta}_i = \omega_i + k\sum_{j=1}^{n} A_{ij} \sin(\theta_j - \theta_i), \quad i = 1,\ldots,n. \tag{1.1}$$

where $\theta_i$ are the phase angles of the $n$ oscillating nodes, $\omega_i$ are the natural frequencies of the nodes and $k$ is the coupling constant. Note that each $\theta_i$ is understood to denote a function of $t$. Where we wish to denote $\theta_i$ at a particular time $t$ we shall make this explicit as $\theta_i(t)$. An important special case is that in which the frequencies $\omega_i$ are all equal. We shall henceforth refer to this case as *homogeneous* and the general case as *inhomogeneous*.

It is well known that results for unique stable fixed points can be obtained for restricted phase angles $\theta_i$ in the range $[-\pi/2, \pi/2]$ (see Jadbabaie, Motee, and Barahona (2004), and Ochab and Gora (2009)). We shall give results for fixed points and their stability where phase angles are over the full range $[-\pi, \pi]$.

We note that by summing over equation (1.1) that the sine terms cancel and the average frequency of the nodes is a constant $\bar{\omega}$ so that (see Strogatz (2000))

$$\frac{1}{n}\sum_{j=1}^{n} \dot{\theta}_j = \frac{1}{n}\sum_{j=1}^{n} \omega_j = \bar{\omega}. \tag{1.2}$$

This result is useful in the following section.

*1.2 Synchronization and fixed points*

It has been observed that for arbitrary initial phases (the $\theta_i(0)$) networks *synchronize* in that some of the node phases converge to the same, or nearly the same phase angle, while the frequencies converge to a common value. Meanwhile the remaining nodes behave non-uniformly or "drift". Moreover at a sufficiently large critical coupling constant $k$, it becomes possible for all nodes to participate in synchronized behaviour.

Define a *frequency fixed point* as a situation in which all the node frequencies are equal and fixed over time. By equation (1.2) this is characterised by

$$\dot{\theta}_i = \bar{\omega}, \quad i = 1,\ldots,n. \tag{1.3}$$

It follows that for a frequency fixed point all phase differences remain constant. This notion is studied in Jadbabaie, Motee, and Barahona (2004), and Wiley, Strogatz and Girvan (2006) (where they refer to *phase locked* fixed points), and Ochab and Gora (2009) (referring to *stationary* fixed points).

We note that a frequency fixed point may also be accompanied by phase synchrony whereby the nodes have the same, or nearly the same phase angles. This is however not necessarily the case, and the phase angles of a frequency fixed point may be significantly different. Examples of this kind are the ring networks studied by Ochab and Gora (2009). The possibility of frequency fixed points away from phase synchrony, and their stability, is the subject of this paper.

Combining equations (1.1) and (1.3) the phase differences at the frequency fixed point satisfy,

$$\bar{\omega} = \omega_i + k \sum_{j=1}^{n} A_{ij} \sin(\theta_j - \theta_i), \quad i = 1,\ldots,n. \tag{1.4}$$

Henceforth we shall use the term fixed point to refer to a frequency fixed point satisfying equation (1.4). The *zero fixed point* refers to the fixed point in which the phase angle differences are all zero. A *non-zero fixed point* refers to a fixed point in which some phase angle difference is non-zero. The character, number and location of the fixed points of a network are clearly important in understanding the types of dynamics that are possible from all possible initial phases of the system (that is where the $\theta_i$ are considered in the full range $-\pi$ to $+\pi$). In particular, understanding the nature of all stable fixed points is particularly important, since each represents a different attractor set of positive $n$-dimensional volume within the set of all states. Note that fixed points can have a range of attractor set types ranging from single points for unstable fixed points, through $m$-dimensional ($m<n$) saddle point structures for partially stable fixed points, to $n$-dimensional volumes for stable fixed points. In this sense we shed some light on understanding the stability of the Kuramoto model globally, a topic with few if any theoretical results (see the remarks at the end of Section 3 of Strogatz (2000)). In this regard we note the work of Diaz-Guilera and Arenas (2008) which includes a discussion of ring networks with equal natural frequencies, and Ochab and Gora (2009) on ring networks with general natural frequencies. Also Wiley, Strogatz and Girvan (2006) study a family of ring like structures in which the natural frequencies are equal, with some analytic results pertaining to certain classes of fixed points as well as numerical simulations. These results are referred to in more detail in Section 4.

*1.3 Multiple stable fixed points*

We emphasise that the number of stable fixed points that may be associated with some networks is potentially high, and this is demonstrated with a simple construction. That ring networks admit multiple stable fixed points has been known to researchers for some time (referred in passing in Jadbabaie, Motee, and Barahona (2004), and studied in more detail by Wiley, Strogatz and Girvan (2006), Diaz-Guilera and Arenas (2008), and Ochab and Gora (2009)). We begin with such a six point ring network that, according to Ochab and Gora (2009) admits two stable fixed point configurations (Figure 1). In this network the natural frequencies are all equal and the node phases relative to a "grounding" node $x$ are indicated. The ring is drawn so as to position the nodes around the rim of a circle, at angles that match the phase angle differences relative to $x$. These configurations are clearly fixed points, but are also stable by virtue of the phase angle differences across the links lying between $-\pi/2$ and $\pi/2$ (see Ochab and Gora (2009)).

By attaching six point rings to nodes of a tree we can then construct networks with an exponential number of stable fixed points where each ring can independently be in one of two states corresponding to the two states of Figure 1. In Figure 2 we show this construction for a star network with a particular choice of stable phase angles indicated for two of the rings.

In the following section we formally develop and investigate stability conditions in the "mixed" case where some of the phase angle differences are between $-\pi/2$ and $\pi/2$ and some are not.

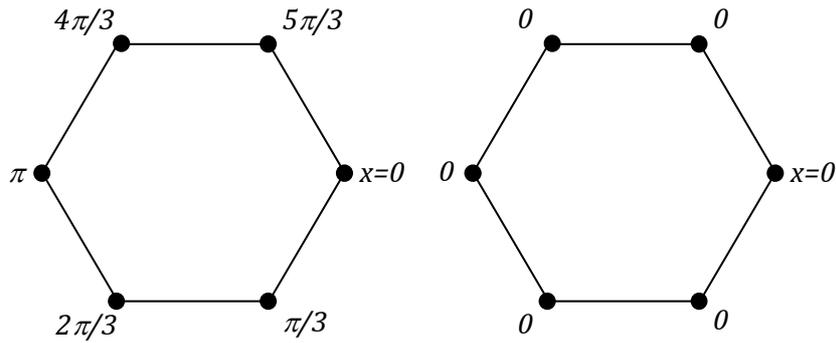

Figure 1 – ring network with two stable fixed points.

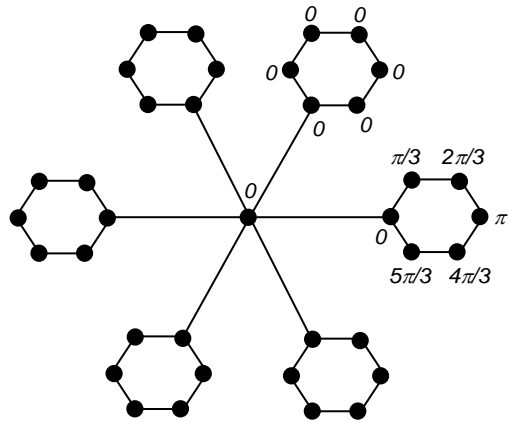

Figure 2 – Network with $2^6$=64 stable fixed points.

*1.4 Roadmap*

We sketch the flow of the paper. The introduction discusses the Kuramoto model, and focuses on the particular part of the problem space addressed here. That is the study of stable fixed points in the homogeneous case of equal natural frequencies with general topology. In Section 2 we develop a necessary condition for the existence of stable fixed points for the general network Kuramoto model. The equations for the homogeneous case are introduced in Section 3, and the significance of the homogeneous case to the general inhomogeneous case is explained. In Section 4 we apply the fixed point stability conditions to show that for the complete network the homogeneous model the only stable fixed point solution over the full space of phase angles in the range $-\pi$ to $\pi$ is the *zero fixed point* in which all the phase angle differences are zero. This result is used to prove the conjecture of Verwoerd and Mason (2007) that for the complete network with homogeneous model there is one fixed point and this has a basin of attraction consisting of the entire space minus a set of measure zero. A more complex argument generalises the stable fixed point result to networks with all node degrees at least $0.9395(n-1)$. The necessary conditions are also tested to see how close to sufficiency they might be by applying them to a class of regular degree networks studied by Wiley, Strogatz and Girvan (2006). This suggests the necessary conditions for stable fixed points are a useful tool in the study of stability. In Section 5 we summarise our results and highlight some outstanding challenges.

## 2 Stability

We develop conditions for the linear stability of fixed points of the Kuramoto model. At this point it is useful to make functions $f$ of $t$ explicit by using $f(t)$. Let $\{\theta_i^*, i=1,...,n\}$ be a fixed point so that by equation (1.3)

$$\dot{\theta}_i^*(t) = \overline{\omega}, \quad i=1,...,n. \tag{2.1}$$

Since $\theta_j^*(t) - \theta_i^*(t)$ is constant for all $i,j$ we shall use the compact notation $\theta_j^*(t) - \theta_i^*(t) = \Delta^*_{ji}$. Near the fixed point we express $\theta_i$ as

$$\theta_i(t) = \theta_i^*(t) + \sigma_i(t)$$

where $\sigma_i(t)$ is small. We then have

$$\theta_j(t) - \theta_i(t) = \theta_j^*(t) - \theta_i^*(t) + \sigma_j(t) - \sigma_i(t) = \Delta^*_{ji} + \sigma_j(t) - \sigma_i(t).$$

By using the first order approximation for sine

$$\sin(\theta_j(t) - \theta_i(t)) \approx \sin(\Delta^*_{ji}) + \cos(\Delta^*_{ji})(\sigma_j(t) - \sigma_i(t)).$$

Substituting into equation (1.1) and using equation (2.1) gives the first order approximation

$$\dot{\sigma}_i(t) \approx -\overline{\omega} + \omega_i + k\sum_{j=1}^{n} A_{ij}[\sin(\Delta^*_{ji}) + \cos(\Delta^*_{ji})(\sigma_j(t) - \sigma_i(t))],$$

$$= -\overline{\omega} + \omega_i + k\sum_{j=1}^{n} A_{ij}\sin(\Delta^*_{ji}) + k\sum_{j=1}^{n} A_{ij}[\cos(\Delta^*_{ji})(\sigma_j(t) - \sigma_i(t))],$$

$$= c_i + k\sum_{j=1}^{n} A_{ij}\cos(\Delta^*_{ji})\sigma_j(t) - k\left[\sum_{j=1}^{n} A_{ij}\cos(\Delta^*_{ji})\right]\sigma_i(t),$$

for some constants $c_i$.

Thus the first order linear approximation takes the matrix form

$$\dot{\boldsymbol{\sigma}} \approx \boldsymbol{c} + kM\boldsymbol{\sigma}^T. \tag{2.2}$$

Where $\boldsymbol{c}$ is a constant vector and $M$ is the matrix defined by

$$M_{ij} = A_{ij}\cos(\Delta^*_{ji}) - \delta_{ij}\sum_{k=1}^{n} A_{ik}\cos(\Delta^*_{ki}),$$

where the Kronecker delta function $\delta_{ij} = 1, 0$ for $i = j, i \neq j$ respectively.

The first order differential system (2.2) is stable if and only if all the eigenvalues of the matrix $M$ are less than or equal to zero. The relationship between the network topology $A_{ij}$, the geometry defined by the $\theta_i^*$, and the negative semi-definiteness of $M$ is a highly complex one. However we can deduce a useful necessary condition for the stability of a fixed point. We note that this condition is not sufficient and so not equivalent to fixed point stability, but follows from it. We explore how close to sufficiency this condition might be in Section 4.3.

**Lemma 2.1** Let $\{\theta_i^*, i=1,\ldots,n\}$ be any stable fixed point solution to the Kuramoto model. Then we have the equalities,

$$\bar{\omega} = \omega_i + k\sum_{j=1}^{n} A_{ij} \sin(\Delta_{ji}^*), \quad i=1,\ldots,n, \tag{2.3}$$

while for any non-empty node subset $S$,

$$\sum_{(i,j)\in(S,S^c)} A_{ij} \cos(\Delta_{ji}^*) \geq 0. \tag{2.4}$$

*Remark.* It follows trivially from this lemma that if for every link ($A_{ij}=1$), $|\Delta^*_{ji}|>\pi/2$ then the fixed point $\{\theta_i^*, i=1,\ldots,n\}$ is unstable. It is also well known (Ochab and Gora (2009)) that if for each link $|\Delta^*_{ji}|<\pi/2$ the fixed point is stable, while if there are phase differences with absolute values both less than $\pi/2$ and greater than $\pi/2$ then the fixed point may be stable or unstable. Lemma 2.1 provides a tool to help resolve this mixed case. The hope is that this condition can more readily be related to the network topology $A_{ij}$ than the negative semi-definiteness of $M$. We emphasise that this condition is a necessary condition only. However in Section 4.3 we present evidence that indicates how close to sufficiency this lemma might be.

**Proof.** Equations (2.3) follow from equation (1.4). Since the system $\{\theta_i^*, i=1,\ldots,n\}$ is stable it follows that the matrix $M$ is negative semi-definite. Thus $z^T M z \leq 0$ for any non-zero vector $z$, or equivalently $\max\{z^T M z\} \leq 0$. Choose $z$ to be the vector defined by $z_i=1$ for $i \in S$ while $z_i=0$ for $i \notin S$. Observe that $z^T M z$ is the sum of elements of $M$ where both the row and column indices are elements of $S$. Thus

$$0 \geq z^T M z = \sum_{i\in S}\sum_{j\in S} M_{ij}$$

$$= \sum_{i\in S}\sum_{j\in S}\left[A_{ij}\cos(\Delta_{ji}^*) - \delta_{ij}\sum_{k=1}^{n} A_{ik}\cos(\Delta_{ji}^*)\right]$$

$$= \sum_{i\in S}\left[\sum_{j\in S} A_{ij}\cos(\Delta_{ji}^*) - \sum_{j\in S}\delta_{ij}\sum_{k=1}^{n} A_{ik}\cos(\Delta_{ki}^*)\right]$$

$$= \sum_{i\in S}\left[\sum_{j\in S} A_{ij}\cos(\Delta_{ji}^*) - \sum_{k=1}^{n} A_{ik}\cos(\Delta_{ji}^*)\right]$$

$$= -\sum_{i\in S}\left[\sum_{k=1}^{n} A_{ik}(\Delta_{ji}^*) - \sum_{j\in S} A_{ij}\cos(\Delta_{ji}^*)\right]$$

$$= -\sum_{i\in S}\sum_{j\in S^c} A_{ij}\cos(\Delta_{ji}^*)$$

$$= -\sum_{(i,j)\in(S,S^c)} A_{ij}\cos(\Delta_{ji}^*).$$

This proves the Lemma. □

**3 Homogeneous systems**

For a fixed point system $\{\theta_i^*, i=1,\ldots,n\}$ we can rearrange equation (2.3) as

$$\frac{\bar{\omega} - \omega_i}{k} = \sum_{j=1}^{n} A_{ij} \sin(\Delta^*_{ji}), \quad i = 1,\ldots,n. \tag{3.1}$$

It follows that as the coupling constant increases in relation to the natural frequencies the left hand side tends to 0, in the limit giving the equation (3.2),

$$0 = \sum_{j=1}^{n} A_{ij} \sin(\Delta^*_{ji}), \quad i = 1,\ldots,n. \tag{3.2}$$

Notice that this limiting equation is the same as that obtained in the homogeneous case, and also that this equation is independent of the value of the equal frequency. Thus as the coupling constant increases, the fixed points for unequal natural frequencies given by (3.1) converge to the solutions for the homogeneous system of equation (3.2). Furthermore in the homogeneous case it is sufficient to consider the case where the natural frequencies are all zero. For this reason the fixed point character of the case with all natural frequencies zero is fundamental to understanding the limiting behaviour of the general inhomogeneous case. The homogeneous case can in some sense be considered a gauge of the contribution of the purely topological effects on fixed points and stability.

Though equation (3.2) can have many (indeed infinitely many) fixed points, we show for a class of dense networks there is a single stable fixed point solution, namely the zero fixed point given by $\Delta^*_{ji} = 0$ for all $i$ and $j$.

In the following section we apply Lemma 2.1 to show that for the complete network there are no non-zero stable fixed point solutions to any homogeneous Kuramoto model. In other words the only stable fixed point solution is the equal phase or zero solution. Using a more complex argument we extend this result to a class of dense networks. As we have shown for the homogeneous system it is sufficient to consider the system where all the natural frequencies are zero. From herein we make this assumption.

## 4 Stability results

*4.1 Complete homogeneous networks*

It is useful to visualise the phase angles as points on the unit circle within the complex plane where the nodes are positioned around a circle in such a way that the phase angles associated with the nodes match their position around the circle. Angles are measured relative to a reference direction chosen to be due east of the centre, with positive angles in a clockwise direction. Thus phase angle differences are equal to the angle subtended at the centre. We illustrate this visualisation in Figure 3 and use it to discuss the utility of the stability Lemma 2.1. In this figure we have circle diagrams for fixed points associated with the hexagonal ring, the complete network on 3 nodes, and a 3-regular 8 node network (in which each node is adjacent to three others). The hexagonal ring network is stable since the phase angle differences are π/3 (Ochab and Gora (2009)). We shall use Lemma 2.1 to show that the 3 and 8 node networks are unstable.

For the 3 node network we set $S=\{a\}$. The left hand side of inequality (2.4) is

$$2\cos(2\pi/3) = -1,$$

and the associated fixed point is unstable by Lemma 2.1. For the 8 node network we take $S$ to be the 4 nodes on one side of the dotted line of Figure 3 (say at π/2, 3π/4, π, and 5π/4). The left hand side of inequality (2.4) is then

$$2\cos(\pi/4) + 4\cos(\pi) = \sqrt{2} - 4 < 0,$$

and again the associated fixed point is unstable by Lemma 2.1.

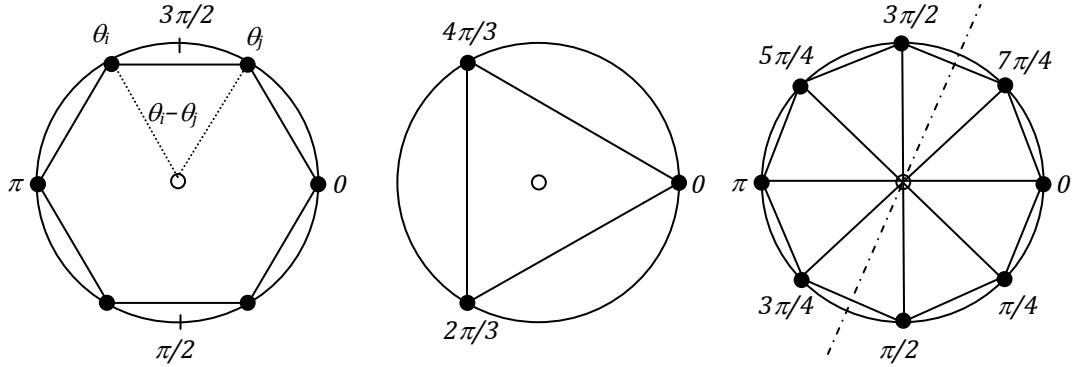

Figure 3- Fixed point networks on 6, 3, and 8 nodes.

It is well known that the zero fixed point is stable, see Arenas *et al* (2008). We show that for any homogeneous system this is the only stable fixed point for the complete network. We emphasise that this is for all phase angles in the full range $[-\pi, \pi]$ rather than the restricted range $[-\pi/2, \pi/2]$. We later extend this result to dense networks. Throughout the rest of this paper we shall use $|A|$ to represent the number of elements in a node set $A$ and $\|x\|$ to represent the modulus of a complex number $x$.

Formalising the representation of nodes on the complex plane let $z_i$, $i=1,...,n$ be given by $z_j = e^{i\theta_j^*} = \cos(\theta_j^*) + i\sin(\theta_j^*)$ (see Figure 4).

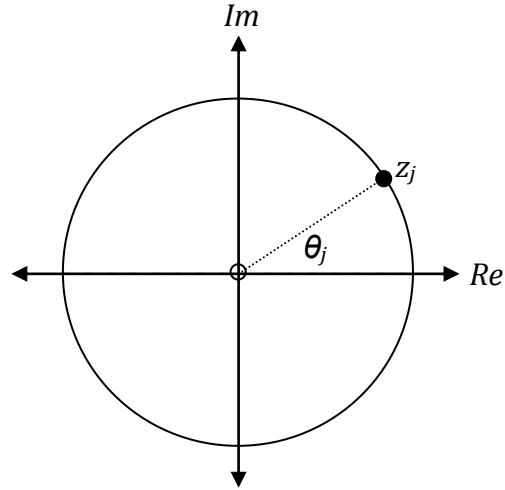

Figure 4

**Theorem 4.1** Any complete network homogeneous Kuramoto model has no non-zero stable fixed point.

**Proof.** Let $\{\theta_i^*, i=1,...,n\}$ be any non-zero stable fixed point solution. Let

$$p = \sum_{i=1}^{n} z_i = \sum_i \cos(\theta_i^*) + i\sum_i \sin(\theta_i^*).$$

Then for any $m$

$$\bar{z}_m p = \bar{z}_m \sum_{i=1}^{n} z_i = \sum_{i=1}^{n} \bar{z}_m z_i = 1 + \sum_{i \neq m} \bar{z}_m z_i$$
$$= 1 + \sum_{i \neq m} \cos(\Delta^*_{im}) + i \sum_{i} \sin(\Delta^*_{im}) = 1 + \sum_{i \neq m} \cos(\Delta^*_{im}),$$
(4.1)

since the imaginary term disappears due to the fixed point condition (3.2).

We consider two cases.

Case $p \neq 0$. For some angle $\theta^p$,

$$p = \|p\|[\cos\theta^p + i\sin\theta^p].$$

Thus

$$\bar{z}_m p = \|p\|[\cos(\theta^p - \theta^*_m) + i\sin(\theta^p - \theta^*_m)].$$
(4.2)

By comparing equations (4.1) and (4.2) we must have $\sin(\theta^p - \theta_m^*) = 0$ for each $m$. It follows that $\theta_m^* = \theta^p$ or $\theta^p + \pi$ for every $m$. If $\theta_m^* = \theta^p$ for all $m$ or $\theta_m^* = \theta^p + \pi$ for all $m$ then we have the zero fixed point solution. If the set of node angles $\theta_m^*$ are divided between $\theta^p$ and $\theta^p + \pi$ with at least one angle of each type let $A$ denote the set of nodes with corresponding angles $\theta^p$, so that $A^c$ is the set of nodes with phase angles $\theta^p + \pi$. We then have non empty node sets $A$ and $A^c$ with

$$\sum_{(i,j) \in (A, A^c)} A_{ij} \cos(\Delta^*_{ji}) = -|A||A^c| \leq -1.$$

This is unstable by inequality (2.4).

Case $p=0$. In this case by equation (4.1) for every $m$

$$\sum_{i \neq m} \cos(\Delta^*_{im}) = -1.$$

Again this case is unstable by equation (2.3), which completes the proof. □

Wiley, Strogatz and Girvan (2006) use a gradient method to show that for a certain class of regular networks (see Section 4.3), including the complete network, attractors must take the form of fixed points. In other words attractors in the form of limit cycles, limit tori or other exotic structures are not possible (see Lorenz for a discussion of these phenomena (1996)). This in turn means that for the complete network the basin off attraction of stable fixed points account for the entire space minus a set of measure zero (the set of measure zero taking the form of partially stable $m<n$ dimensional saddle point structures and (totally) unstable fixed points). By Theorem 4.1 the homogeneous complete network only has the zero stable fixed. This proves a conjecture of Verwoerd and Mason (2007).

**Corollary 4.1** For the complete network the zero fixed point solution to any corresponding homogeneous Kuramoto model has a basin of attraction consisting of the entire space minus a set of measure zero.

*4.2 Dense Networks*

We generalise Theorem 4.1 to include all networks with node degrees at least $\mu(n-1)$ where $\mu$ is a constant less than 1, $n-1$ being the maximum possible degree. We suspect that this result is far from the best possible of this kind (the theorem gives our best lower bound of 0.9395 on $\mu$). The author has been unable to find a better bound to what appears to be a difficult geometrical problem. The proof is by contradiction. The constant 0.9395 is chosen to be as small as possible while allowing the contradiction to follow (in which the right hand side of inequality (4.8) is negative). This constant also depends upon the transcendental solution to the equation $\tan(x)=1/x$ (see Case 2 of the proof), and so we would not expect it to have a closed form or surd representation.

**Theorem 4.2** Consider any $n$ node network with node degrees at least $0.9395(n-1)$. Then any corresponding homogeneous Kuramoto model has no non-zero stable fixed point.

**Proof.** Let $\{\theta_i^*, i=1,\ldots,n\}$ be any non-zero stable fixed point solution. We shall show that this leads to a contradiction. Using complex arithmetic we have $z_j = e^{i\theta_j^*}$, and $p = \sum_{i=1}^{n} z_i$. Recall by equations (4.1) and (4.2) that for any $m$

$$\bar{z}_m p = \sum_i \cos(\Delta_{im}^*) + i \sum_i \sin(\Delta_{im}^*) = \|p\|[\cos(\theta^p - \theta_m^*) + i\sin(\theta^p - \theta_m^*)].$$

Now the node $m$ is adjacent to at least $0.9395(n-1)$ of the nodes of $A^c$, and not adjacent to at most $0.0605(n-1)$ nodes of $A^c$. Using the fixed point equation (3.2) it follows that

$$\left|\sum_i \sin(\Delta_{im}^*)\right| \leq \left|\sum_i A_{mi}\sin(\Delta_{im}^*)\right| + 0.0605(n-1) = 0.0605(n-1).$$

Thus,

$$\|p\|\left|\sin(\theta_m^* - \theta^p)\right| \leq 0.0605(n-1). \tag{4.3}$$

We consider two cases.

Case 1 $\|p\| > \sqrt{2} \times 0.0605(n-1)$. By inequality (4.6) $\sin(\theta^p - \theta_m^*) < 1/\sqrt{2}$. This in turn means that

$$-\pi/4 < \theta_m^* - \theta^p < \pi/4 \text{ or } -5\pi/4 < \theta_m^* - \theta^p < 3\pi/4.$$

This means that all nodes are contained within two node sets $A$ and $B$ associated with arcs of less than $\pi/2$ as shown in Figure 5. Let us assume that one of $A$ or $B$ is empty, say $B$. Let $g$ be any node closest to the edge of the arc $A$. Consider equality (3.2) where $i$ corresponds to the node $g$. Then each of the terms in the sum of (3.2) has the same sign, and the terms are not all zero unless we have the zero solution. Thus (3.2) cannot be satisfied at $i=g$ and we cannot have a fixed point.

Assume on the other hand that both $A$ and $B$ are not empty. Now by an elementary result of network theory any network with node degrees at least $\frac{1}{2}(n-1)$ must be connected (see Bondy and Murty (1976)). It follows that there must be at least one edge between $A$ and $B$. The cosine of the phase difference of all edges between *A* and *B* is clearly negative so that inequality (2.4) is violated. The system is therefore unstable.

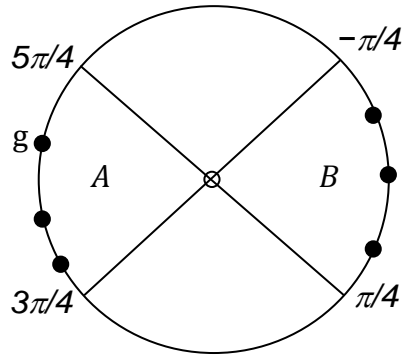

Figure 5

Case 2 $\|p\| \leq \sqrt{2} \times 0.0605(n-1) < 0.0856(n-1)$. Let *A* be the nodes contained in an arc of $r=1.7206$ radians containing a maximum number of nodes (more precisely $r=2x$ where $x$ is the solution to $\tan(x)=1/x$). Then by an averaging argument *A* must contain at least $rn/2\pi=0.2738n$ nodes. Let *v* be the unit vector in the complex plane that bisects *A*. Then the angles between *v* and any $z_i$ in *A* are between plus or minus $r/2=0.8603$ radians (see Figure 6).

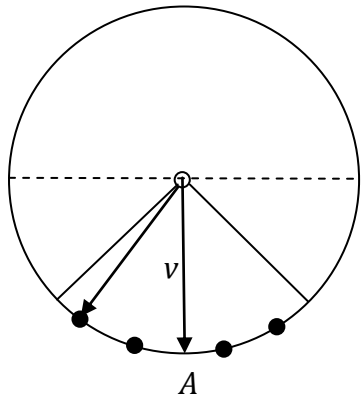

Figure 6

Let $p_A = \sum_{i \in A} z_i$ and set $|A|=\beta n$ where, as we have shown, $\beta \geq 0.2738$. Then

$$\|p_A\| \geq v \bullet p_A = \sum_{i \in A} v \bullet z_i \geq \cos(0.8603) \, |A| = 0.6522\beta n > 0.1785n. \tag{4.4}$$

Now

$$\text{Re}[\bar{p}_A p] = \text{Re}\left[\sum_{i \in A} \bar{z}_i \sum_{i=1}^{n} z_i\right] = \text{Re}\left[\sum_{i \in A} \bar{z}_i \left(\sum_{i \in A} z_i + \sum_{i \in A^c} z_i\right)\right]$$

$$= \text{Re}\left[\sum_{i \in A} \bar{z}_i \sum_{i \in A} z_i\right] + \text{Re}\left[\sum_{i \in A} \bar{z}_i \sum_{i \in A^c} z_i\right] \quad (4.5)$$

$$= \left\|\sum_{i \in A} z_i\right\|^2 + \sum_{(i,j) \in (A,A^c)} \cos(\Delta_{ij}^*)$$

$$= \|p_A\|^2 + \sum_{(i,j) \in (A,A^c)} \cos(\Delta_{ij}^*).$$

Each node of $A$ is adjacent to at least $0.9395(n-1)$ of the nodes of $A^c$, and not adjacent to at most $0.0605(n-1)$ nodes of $A^c$. It follows that

$$\sum_{(i,j) \in (A,A^c)} A_{ij} \cos(\Delta_{ij}^*) \le \left[\sum_{(i,j) \in (A,A^c)} \cos(\Delta_{ij}^*)\right] + 0.0605(n-1)|A|. \quad (4.6)$$

Combining equation (4.5) with inequality (4.6) we have

$$\sum_{(i,j) \in (A,A^c)} A_{ij} \cos(\Delta_{ij}^*) \le \text{Re}[p_A p] - \|p_A\|^2 + 0.0605(n-1)|A|$$

$$\le \|p_A\| \times \|p\| - \|p_A\|^2 + 0.0605(n-1)|A| \quad (4.7)$$

$$\le \|p_A\|(\|p\| - \|p_A\|) + 0.0605(n-1)|A|.$$

Now by the inequality defining this case and the bound (4.4)

$$\|p\| - \|p_A\| \le 0.0856(n-1) - 0.1785n < 0.$$

It follows that we can substitute the upper bounds for $\|p\|$ and the lower bounds for $\|p_A\|$ into inequality (4.7) to obtain

$$\sum_{(i,j) \in (A,A^c)} A_{ij} \cos(\Delta_{ij}^*) \le \|p_A\|(\|p\| - \|p_A\|) + 0.0605(n-1)|A|$$

$$\le 0.6522\beta n(0.0856(n-1) - 0.6522\beta n) + 0.0605\beta n(n-1)$$

$$\le 0.6522\beta n(0.0856n - 0.6522\beta n) + 0.0605\beta n^2$$

$$\le 0.2738n^2(0.0856 \times 0.6522 - 0.6522^2 \times 0.2738 + 0.0605)$$

$$\le -0.00002n^2. \quad (4.8)$$

This last inequality contradicts the stability inequality (2.4) which completes the proof. □

*4.3 The Wiley Strogatz Girvan Networks*

It is interesting to understand the limitations of Lemma 2.1 as a necessary but not necessarily sufficient condition for stability. To do this we draw on the analysis and of Wiley, Strogatz and Girvan (2006) as applied to a particular network class. Consider the homogeneous regular $d$-degree network on $n$ nodes in which the nodes are adjacent to the nearest $d$ neighbours around a circle ($d$ even). We shall call these *WSG* networks. Furthermore we make the plausible simplifying assumption (Wiley, Strogatz and Girvan (2006)) that stable fixed points for such networks are characterised by equal phase differences between adjacent nodes, called "uniformly twisted" states (see Figure 7).

Assume then that we have a non-zero stable fixed point which is uniformly twisted. Inequality (2.4) is applied. Consider the node partition $A$ formed by dividing the network in half (see Figure 7). The summation term of inequality (2.4) involves the sum of phase difference terms between nodes of $A$ and $A^c$. There are 2 such terms with phase difference $\delta=2\pi/n$, 4 terms with phase differences $2\delta$, 6 terms with phase difference $3\delta$,.., up to $d$ terms with phase differences $d\delta/2$ (see Figure 7). We therefore have

$$\sum_{(i,j)\in(A,A^c)} A_{ij}\cos(\Delta^*_{ij}) = 2\sum_{i=0}^{d/2} i\cos(i\delta), \text{ where } \delta = \frac{2\pi}{n}.$$

For large $n$ the sum approaches the corresponding integral so that,

$$2\sum_{i=0}^{d/2} i\cos(i\delta) = \frac{2}{\delta^2}\sum_{i=0}^{d/2} i\delta\cos(i\delta)\delta \to \frac{n^2}{2\pi^2}\int_0^{d\pi/n} x\cos(x)dx = \frac{n^2}{2\pi^2}\left[\frac{d\pi}{n}\sin(d\pi/n) + \cos(d\pi/n) - 1\right].$$

Thus the limiting form of inequality (2.4) becomes

$$\frac{n^2}{2\pi^2}\left[\frac{d\pi}{n}\sin(d\pi/n) + \cos(d\pi/n) - 1\right] > 0.$$

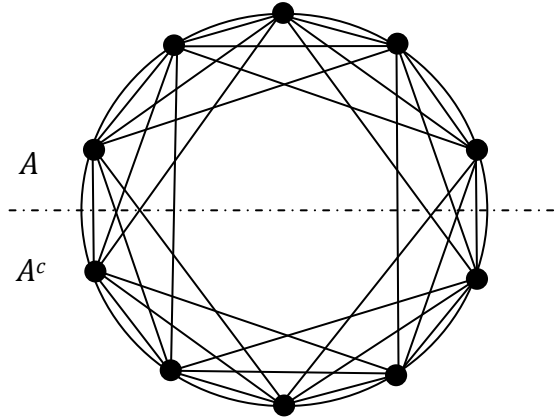

Figure 7 – regular degree network

This has the solution $d\pi/n<2.3311$ or $d/n<0.7420$. This means that the limiting maximum value of $d/n$ allowing stable uniformly twisted fixed points is less than 0.7420. Wiley, Strogatz and Girvan (2006) show this limiting value to be 0.6809 (see Table 1 of their paper in which the class of stable fixed points we describe here correspond to those with a single full twist in the state, that is the $q=1$ case). We have also verified this by computing the eigenvalues of the corresponding matrix $M$ from equation (2.2) for large $n$ (=100) and varying $d$. This suggests that the necessary conditions (2.3) and (2.4) have the potential to be useful for the study of stable fixed points.

It seems that a definitive topological necessary and sufficient condition for the existence of a non-zero stable fixed point is a complex problem even in the case of the homogeneous system.

At this point we note that the gradient result of Wiley, Strogatz and Girvan (2006) combined with Theorem 4.2 gives the analogous result to Corollary 4.1.

**Corollary 4.2** For any *WSG* network with node degrees at least $0.9395(n-1)$ the zero fixed point solution to any corresponding homogeneous Kuramoto model has a basin of attraction consisting of the entire space minus a set of measure zero.

## 5 Conclusion

Characterising the stable fixed points of the Kuramoto model is central to understanding the global dynamics of the model, since each stable fixed point has attractor sets of positive measure over the complete phase space. Furthermore the homogeneous Kuramoto model, apart from being of interest in its own right, represents the behaviour of the inhomogeneous model in the limit as the coupling constant tends to infinity. Thus the characterization of the stable fixed points in the homogeneous case is a valuable place to start in understanding the general (inhomogeneous) case. Despite this, and the interest in understanding this problem (see Jadbabaie, Motee, and Barahona (2004), Wiley, Strogatz and Girvan (2006), Ochab and Gora (2009), and Verwoerd and Mason (2007)), there is a dearth of results concerning stable fixed points over the full space of phase angles $[-\pi,\pi]$, specifically the conditions on network topology for non-zero stable fixed points.

We have presented a new necessary condition for stable fixed points, and used this to show that for the homogeneous system there is only one stable fixed point solution for all networks with sufficiently high node degrees, namely the zero solution (see Theorems 4.1 and 4.2). This result proves that for a class of dense regular networks the zero fixed point of the homogeneous Kuramoto model has an attractor basin consisting of the entire space minus a set of measure zero.

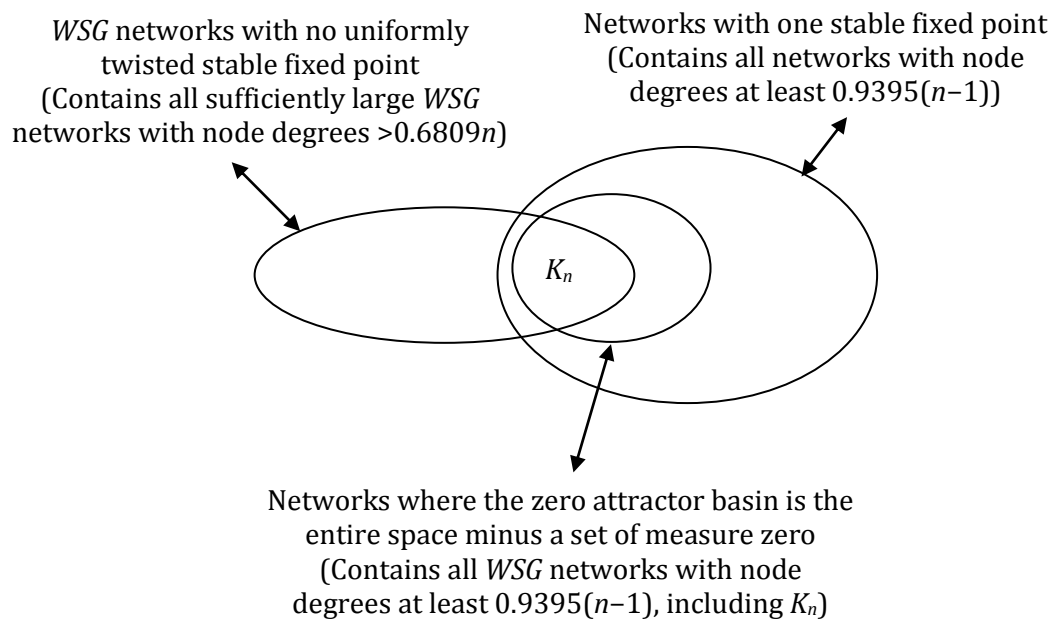

Figure 8 – Venn diagram summary of results

Many of the results in this paper are summarised in the Venn diagram of Figure 8 showing the relationships between the networks with one stable fixed point, the *WSG* networks, and the networks where the zero attractor basin is the entire space minus a set of measure zero. $K_n$ is the complete network on *n* nodes.

Much remains to be understood however about the relationship between network topology, and the associated fixed points or other potential attractor sets in both the homogeneous and inhomogeneous Kuramoto model, as well understanding as the geometry of the corresponding basins of attraction. Is it possible, for example, that attractors in the form of limit cycles, limit tori or other exotic structures are possible for the homogeneous model? In the context of improving the main result of this paper, what is the smallest constant $\mu$ such that all networks with node degrees at least $\mu(n-1)$ have a single stable fixed point? Our results together with those of Wiley, Strogatz and Girvan (2006) show that $0.6809 \leq \mu \leq 0.9395$.

## Acknowledgements

The author is thankful to colleagues with whom a number of fruitful conversations have taken place over the course of this research. They are A. Dekker, M. Sweeney and A. Kalloniatis, the latter also providing detailed feedback on an early draft of this work. Thanks also to one of the referees who made a number of recommendations that led to improvements in the paper.